# ASYMPTOTIC THEORY OF SEMIPARAMETRIC $Z$-ESTIMATORS FOR STOCHASTIC PROCESSES WITH APPLICATIONS TO ERGODIC DIFFUSIONS AND TIME SERIES

By Yoichi Nishiyama

*Institute of Statistical Mathematics*

This paper generalizes a part of the theory of $Z$-estimation which has been developed mainly in the context of modern empirical processes to the case of stochastic processes, typically, semimartingales. We present a general theorem to derive the asymptotic behavior of the solution to an estimating equation $\theta \rightsquigarrow \Psi_n(\theta, \widehat{h}_n) = 0$ with an abstract nuisance parameter $h$ when the compensator of $\Psi_n$ is random. As its application, we consider the estimation problem in an ergodic diffusion process model where the drift coefficient contains an unknown, finite-dimensional parameter $\theta$ and the diffusion coefficient is indexed by a nuisance parameter $h$ from an infinite-dimensional space. An example for the nuisance parameter space is a class of smooth functions. We establish the asymptotic normality and efficiency of a $Z$-estimator for the drift coefficient. As another application, we present a similar result also in an ergodic time series model.

**1. Introduction.** Let us begin with stating our motivating example; the details are presented in Section 4. Consider the one-dimensional ergodic diffusion process $X$ on $I = (l, r) \subseteq \mathbb{R}$ which is a solution to the stochastic differential equation (SDE) given by

$$(1) \qquad X_t = X_0 + \int_0^t S(X_s; \theta) \, ds + \int_0^t \sigma(X_s; h) \, dW_s,$$

where $s \rightsquigarrow W_s$ is a standard Brownian motion. Here, we consider a $d$-dimensional parametric family $\{S(\cdot; \theta); \theta \in \Theta\}$ for the drift coefficient indexed by a compact subset $\Theta$ of $\mathbb{R}^d$, and a possibly infinite-dimensional "parametric" family $\{\sigma^2(\cdot; h); h \in H\}$ for the diffusion coefficient indexed by a (general) totally bounded metric space $(H, d_H)$. We denote by $(\theta_0, h_0)$ the true value of $(\theta, h)$. Our aim is to estimate $\theta_0$ when the model is perturbed by the unknown nuisance parameter $h$. As for the parameter $h_0$, we construct a $d_H$-consistent estimator $\widehat{h}_n$. We prove that the $Z$-estimator $\widehat{\theta}_n$, which is a solution to an estimating equation $\Psi_n(\theta, \widehat{h}_n) = 0$, is asymptotically normal and efficient. [We follow van der Vaart and Wellner (1996) for the terminology "$Z$-estimator."]

There exist a lot of works which treat the estimation problem for the drift coefficient. It is well known that when the process $X = (X_t)_{t \in [0,\infty)}$ is observed *continuously* on the time interval $[0, T]$, the diffusion coefficient may be assumed to be known without loss of generality. (So, we may put $h = h_0$.) In such cases, the asymptotic normality and efficiency of the maximum likelihood estimator (MLE) $\widehat{\theta}_T$ for $\theta_0$, as $T \to \infty$, has been already established. See, for example, Kutoyants (2004). The MLE $\widehat{\theta}_T$ is a solution to the estimating equation $\dot{\ell}_T(\theta) = 0$ with

$$\dot{\ell}_T(\theta) = \frac{1}{T} \int_0^T \frac{\dot{S}(X_t; \theta)}{\sigma^2(X_t; h_0)} [dX_t - S(X_t; \theta) \, dt],$$

where $\dot{S}$ denotes the derivative of $S$ with respect to $\theta$. On the other hand, when the process $X = (X_t)_{t \in [0,\infty)}$ is observed only at discrete time points $\{0 = t_0^n < t_1^n < \cdots < t_n^n\}$, the diffusion coefficient has to be estimated, too. Florens-Zmirou (1989), Yoshida (1992) and Kessler (1997), among others, considered such situations when $H$ is a *finite-dimensional* parameter space, and proved the asymptotic efficiency of some estimators $\widehat{\theta}_n$ for $\theta_0$. Our result does not include these works as special cases, because we assume a condition, which is theoretically strong but practically reasonable, that

$$\Delta_n = \max_{1 \leq i \leq n} |t_i^n - t_{i-1}^n| = o((t_n^n)^{-1}) \quad \text{and} \quad t_n^n \to \infty,$$

which is almost the same as the assumption $n\Delta_n^2 \to 0$. For example, Kessler's (1997) assumption $n\Delta_n^p \to 0$ for a given $p \geq 2$ is weaker than ours. Another difference is that the preceding works derive not only the consistency of the finite-dimensional estimator $\widehat{h}_n$ but also its asymptotic distribution, while we prove only the $d_H$-consistency. However, our work is the first attempt to propose an asymptotically efficient estimator $\widehat{\theta}_n$ for $\theta_0$ when the nuisance parameter $h$ belongs to an *infinite-dimensional* space $(H, d_H)$. Here, by "asymptotically efficient" we mean that the rescaled residual $\sqrt{t_n^n}(\widehat{\theta}_n - \theta_0)$ has the same asymptotic distribution as the continuous observation case, with $h = h_0$ being known, which has been shown to be optimal in the framework of *local asymptotic normality* theory.

We approach this problem by using the approximation of $\dot{\ell}_T(\theta)$ given by

$$\Psi_n(\theta, h) = \frac{1}{t_n^n} \sum_{i=1}^n \frac{\dot{S}(X_{t_{i-1}^n}; \theta)}{\sigma^2(X_{t_{i-1}^n}; h)} [X_{t_i^n} - X_{t_{i-1}^n} - S(X_{t_{i-1}^n}; \theta) |t_i^n - t_{i-1}^n|],$$

where $h_0$ has been replaced by the unknown parameter $h$. Its "compensator" is

$$\widetilde{\Psi}_n(\theta, h) = \frac{1}{t_n^n} \sum_{i=1}^{n} \frac{\dot{S}(X_{t_{i-1}^n}; \theta)}{\sigma^2(X_{t_{i-1}^n}; h)} \left[ \int_{t_{i-1}^n}^{t_i^n} S(X_t; \theta_0) \, dt - S(X_{t_{i-1}^n}; \theta) |t_i^n - t_{i-1}^n| \right],$$

in the sense that the difference $\Psi_n(\theta, h) - \widetilde{\Psi}_n(\theta, h)$ is the terminal variable of a martingale. The key points are to show the weak convergence of the rescaled random fields $(\theta, h) \rightsquigarrow r_n(\Psi_n(\theta, h) - \widetilde{\Psi}_n(\theta, h))$ for some constant $r_n$ tending to $\infty$, and to show the differentiability of $(\theta, h) \rightsquigarrow \widetilde{\Psi}_n(\theta, h)$ around $(\theta_0, h_0)$. Roughly speaking, our main result asserts that if we assume $h \mapsto \sigma^2(\cdot; h)$ is Lipschitz continuous with respect to $d_H$, that the metric entropy condition is satisfied,

$$\int_0^1 \sqrt{\log N(H, d_H, \varepsilon)} \, d\varepsilon < \infty,$$

where $N(H, d_H, \varepsilon)$ is $\varepsilon$-covering number of $H$ with respect to $d_H$, and that we have a $d_H$-consistent estimator $\widehat{h}_n$ for $h_0$, then we can derive the asymptotic distribution of $r_n(\widehat{\theta}_n - \theta_0)$. The consistency of $\widehat{h}_n$ should be established separately.

This approach is based on a new theory for general $Z$-estimators with infinite-dimensional nuisance parameters presented in Section 2, although its proof is just an adaptation of that of Chapter 3.3 of van der Vaart and Wellner (1996) who considered the case where the compensator $\widetilde{\Psi}_n$ is neither random nor depending on $n$. Hopefully, this extension considerably enlarges the application fields of van der Vaart and Wellner's theory to various stochastic process models. Indeed, we also present a result for time series in Section 5, which is briefly introduced below. Kosorok's (2008) new book does not seem to cover our examples.

In Section 5, we will consider an ergodic time series model of the form

$$X_i = S(X_{i-1}, \ldots, X_{i-p}; \theta) + \sigma(X_{i-1}, \ldots, X_{i-q}; h) w_i,$$

where $E[w_i | \mathcal{F}_{i-1}] = 0$ and $E[w_i^2 | \mathcal{F}_{i-1}] = 1$. Here, $\theta$ is an estimated parameter which belongs to a compact subset of $\mathbb{R}^d$, while $h$ is a nuisance parameter from a totally bounded metric space $(H, d_H)$. In the same way as the diffusion process case, we present a general result to derive the asymptotic normality (and efficiency in some cases) of a $Z$-estimator for $\theta_0$. Although there are vast literatures in time series analysis [see, e.g., Taniguchi and Kakizawa (2000)] apparently, our result is new.

The crucial point of our approach is how to show the weak convergence of the random fields $(\theta, h) \rightsquigarrow r_n(\Psi_n(\theta, h) - \widetilde{\Psi}_n(\theta, h))$. For this purpose, we use the general weak convergence theory for $\ell^\infty$-valued martingales established by Nishiyama (1996, 1997, 1999, 2000a, 2000b, 2007). The theory is a good marriage between the martingale theory which has a long history [see, e.g., Jacod and Shiryaev (1987)] and the modern theory of empirical processes [see, e.g., van der Vaart and Wellner (1996)].

The organization of the paper is as follows. In Section 2, we present a general theory for $Z$-estimation with infinite-dimensional nuisance parameter. In Section 3, we prepare a uniform law of large numbers for random fields with abstract parameter, which is often used in the course of our work. The results for the ergodic diffusion process models are presented in Section 4, while those for the ergodic time series models are given in Section 5.

We refer to van der Vaart and Wellner (1996) for the weak convergence theory in $\ell^\infty(T)$-space, where $\ell^\infty(T)$ is the space of bounded functions on a set $T$. We denote by $C_\rho(T)$ the space of functions on $T$ which are continuous with respect to the metric $\rho$. We equip both spaces with the uniform metric. Given a probability measure $P$, we denote by $P^*$ the corresponding outer probability; see van der Vaart and Wellner (1996) for the stochastic convergence theory which does not assume the measurability. We denote by $\xrightarrow{p}$ and $\xrightarrow{d}$ the convergence in (outer) probability and the weak convergence. The limit notation mean in principle that we take the limit as $n \to \infty$. The Euclidean metric on $\mathbb{R}^d$ is denoted by $\|\cdot\|$.

**2. General theory for semiparametric $Z$-estimation.** Let two sets $\Theta$ and $H$ be given. Let

$$\Psi_n : \Theta \times H \to \mathbb{R}^d \quad \text{and} \quad \widetilde{\Psi}_n : \Theta \times H \to \mathbb{R}^d$$

be random maps. The latter should be a random "compensator" of the former, and in the i.i.d. case it is not random and not depending on $n$. Compare the above setting with that in Chapter 3.3 of van der Vaart and Wellner (1996) where $\widetilde{\Psi}_n \equiv \Psi$.

We present a way to derive the asymptotic behaviour of the estimator $\widehat{\theta}_n$ for the parameter $\theta \in \Theta$ of interest, with help of the estimator $\widehat{h}_n$ for the nuisance parameter $h \in H$, which are solutions to the estimating equation

$$\Psi_n(\widehat{\theta}_n, \widehat{h}_n) \approx 0.$$

Here, the true values $\theta_0 \in \Theta$ and $h_0 \in H$ are supposed to satisfy

$$\widetilde{\Psi}_n(\theta_0, h_0) \approx 0.$$

The following theorem extends a special case of Theorem 3.3.1 of van der Vaart and Wellner (1996). See also Theorem 5.21 of van der Vaart (1998).

THEOREM 2.1. *Let $\Theta$ be a subset of $\mathbb{R}^d$ with the Euclidean metric $\|\cdot\|$. Let $(H, d_H)$ be a semimetric space. Let $\Psi_n : \Theta \times H \to \mathbb{R}^d$ and $\widetilde{\Psi}_n : \Theta \times H \to \mathbb{R}^d$ be random maps defined on a probability space $(\Omega_n, \mathcal{F}_n, P_n)$. (We do not* *assume any measurability.) Suppose that there exist a sequence of constants $r_n \uparrow \infty$, some fixed point $(\theta_0, h_0)$ and an invertible matrix $V_{\theta_0, h_0}$ which satisfy the following* (i) *and* (ii).

(i) *There exists a neighborhood $U \subset \Theta \times H$ of $(\theta_0, h_0)$ such that*

$$r_n(\Psi_n - \widetilde{\Psi}_n) \xrightarrow{d} Z \quad \text{in } \ell^\infty(U),$$

*where almost all paths $(\theta, h) \rightsquigarrow Z(\theta, h)$ are continuous with respect to $\rho = \|\cdot\| \vee d_H$.*

(ii) *For given random sequence $(\widehat{\theta}_n, \widehat{h}_n)$, it holds that*

$$\widetilde{\Psi}_n(\widehat{\theta}_n, \widehat{h}_n) - \widetilde{\Psi}_n(\theta_0, h_0) - V_{\theta_0, h_0}(\widehat{\theta}_n - \theta_0) = o_{P_n^*}(r_n^{-1} + \|\widehat{\theta}_n - \theta_0\|)$$

*and that*

$$\|\widehat{\theta}_n - \theta_0\| \vee d_H(\widehat{h}_n, h_0) = o_{P_n^*}(1), \qquad \Psi_n(\widehat{\theta}_n, \widehat{h}_n) = o_{P_n^*}(r_n^{-1}),$$
$$\widetilde{\Psi}_n(\theta_0, h_0) = o_{P_n^*}(r_n^{-1}).$$

*Then it holds that*

$$r_n(\widehat{\theta}_n - \theta_0) \xrightarrow{d} -V_{\theta_0, h_0}^{-1} Z(\theta_0, h_0).$$

To prove the above theorem, we need the following lemma, which is a slight generalization of Lemma 19.24 of van der Vaart (1998).

LEMMA 2.2. *Let $(T, \rho)$ be a semimetric space. Suppose that $Z_n \xrightarrow{d} Z$ in $\ell^\infty(T)$ and that almost all paths of $Z$ are continuous with respect to $\rho$. If $T$-valued random sequence $\widehat{t}_n$ satisfies $\rho(\widehat{t}_n, t_0) = o_{P_n^*}(1)$ for some nonrandom $t_0 \in T$, then $Z_n(\widehat{t}_n) - Z_n(t_0) = o_{P_n^*}(1)$.*

PROOF. Let us equip the space $\ell^\infty(T) \times T$ with the metric $\|\cdot\|_T + \rho$, where $\|\cdot\|_T$ denotes the uniform metric on $\ell^\infty(T)$. Define the function $g: \ell^\infty(T) \times T \to \mathbb{R}$ by $g(z, t) = z(t) - z(t_0)$. Then for any $z \in C_\rho(T)$ and $t \in T$, the function $g$ is continuous at $(z, t)$. Indeed, if $(z_n, t_n) \to (z, t)$, then $\|z_n - z\|_T \to 0$, and thus $z_n(t_n) = z(t_n) + o(1) \to z(t)$, while $z_n(t_0) \to z(t_0)$ is trivial.

By assumption, we have $(Z_n, \widehat{t}_n) \xrightarrow{d} (Z, t_0)$ in $\ell^\infty(T) \times T$ [see, e.g., Theorem 18.10(v) of van der Vaart (1998)]. Since almost all paths of $Z$ belong to $C_\rho(T)$, by the continuous mapping theorem,

$$Z_n(\widehat{t}_n) - Z_n(t_0) = g(Z_n, \widehat{t}_n) \xrightarrow{d} g(Z, t_0) = Z(t_0) - Z(t_0) = 0.$$

The proof is finished. □

PROOF OF THEOREM 2.1. Applying Lemma 2.2 to the $\ell^\infty(U)$-valued random element $Z_n = r_n(\Psi_n - \widetilde{\Psi}_n)$, we have

$$r_n(\Psi_n(\widehat{\theta}_n, \widehat{h}_n) - \widetilde{\Psi}_n(\widehat{\theta}_n, \widehat{h}_n)) - r_n(\Psi_n(\theta_0, h_0) - \widetilde{\Psi}_n(\theta_0, h_0)) = o_{P_n^*}(1).$$

Since $r_n \Psi_n(\widehat{\theta}_n, \widehat{h}_n) = o_{P_n^*}(1)$, we have

$$r_n(\widetilde{\Psi}_n(\widehat{\theta}_n, \widehat{h}_n) - \widetilde{\Psi}_n(\theta_0, h_0)) = -r_n \Psi_n(\theta_0, h_0) + o_{P_n^*}(1).$$

By the assumption (ii), it holds that

$$(2) \qquad r_n V_{\theta_0, h_0}(\widehat{\theta}_n - \theta_0) = -r_n \Psi_n(\theta_0, h_0) + o_{P_n^*}(1 + r_n \|\widehat{\theta}_n - \theta_0\|).$$

Now, since

$$r_n \|\widehat{\theta}_n - \theta_0\| \leq \|V_{\theta_0, h_0}^{-1}\| r_n \|V_{\theta_0, h_0}(\widehat{\theta}_n - \theta_0)\|$$

$$\leq O_{P_n}(1) + o_{P_n^*}(1 + r_n \|\widehat{\theta}_n - \theta_0\|) \qquad \text{by (2)},$$

it holds that $r_n \|\widehat{\theta}_n - \theta_0\| = O_{P_n^*}(1)$. Inserting this to (2), we have

$$r_n V_{\theta_0, h_0}(\widehat{\theta}_n - \theta_0) = -r_n \Psi_n(\theta_0, h_0) + o_{P_n^*}(1)$$
$$= -r_n(\Psi_n(\theta_0, h_0) - \widetilde{\Psi}_n(\theta_0, h_0)) + o_{P_n^*}(1)$$
$$\xrightarrow{d} -Z(\theta_0, h_0),$$

which implies the conclusion. □

In Theorem 2.1, both "$\|\widehat{\theta}_n - \theta_0\| = o_{P_n^*}(1)$" and "$d_H(\widehat{h}_n, h_0) = o_{P_n^*}(1)$" are assumed. Under some conditions, the former automatically follows from the latter, as it is seen in the following theorem.

THEOREM 2.3. *Let $(\Theta, d_\Theta)$ and $(H, d_H)$ be two semimetric spaces. Let $\Psi_n : \Theta \times H \to \mathbb{R}^d$ be a random map defined on a probability space $(\Omega_n, \mathcal{F}_n, P_n)$. (We do not assume any measurability.) Let $\Psi : \Theta \times H \to \mathbb{R}^d$ be a nonrandom function. Suppose that*

$$\sup_{(\theta, h) \in \Theta \times H} \|\Psi_n(\theta, h) - \Psi(\theta, h)\| = o_{P_n^*}(1).$$

*Suppose also that for some $(\theta_0, h_0) \in \Theta \times H$*

$$\inf_{\theta : d_\Theta(\theta, \theta_0) > \varepsilon} \|\Psi(\theta, h_0)\| > 0 \qquad \forall \varepsilon > 0,$$

*and $h \rightsquigarrow \Psi(\theta, h)$ is continuous at $h_0$ uniformly in $\theta$. Then for any random sequence $(\widehat{\theta}_n, \widehat{h}_n)$ such that $\Psi_n(\widehat{\theta}_n, \widehat{h}_n) = o_{P_n^*}(1)$ and that $d_H(\widehat{h}_n, h_0) = o_{P_n^*}(1)$, it holds that $d_\Theta(\widehat{\theta}_n, \theta_0) = o_{P_n^*}(1)$.*

PROOF. Observe that
$$\|\Psi(\widehat{\theta}_n, \widehat{h}_n)\| \leq \|\Psi(\widehat{\theta}_n, \widehat{h}_n) - \Psi_n(\widehat{\theta}_n, \widehat{h}_n)\| + \|\Psi_n(\widehat{\theta}_n, \widehat{h}_n)\|$$
$$\leq \sup_{\theta,h} \|\Psi(\theta, h) - \Psi_n(\theta, h)\| + \|\Psi_n(\widehat{\theta}_n, \widehat{h}_n)\|$$
$$= o_{P_n^*}(1).$$

Now, for every $\varepsilon > 0$, there exist $\delta, \eta > 0$ such that $\|\Psi(\theta, h)\| > \eta$ for every $\theta$ with $d_\Theta(\theta, \theta_0) > \varepsilon$ and every $h$ with $d_H(h, h_0) < \delta$. Thus, the event $\{d_\Theta(\widehat{\theta}_n, \theta_0) > \varepsilon\}$ is contained in the event $\{\|\Psi(\widehat{\theta}_n, \widehat{h}_n)\| > \eta\} \cup \{d_H(\widehat{h}_n, h_0) \geq \delta\}$. The outer probability of the latter event converges to 0. □

**3. Uniform law of large numbers.** In this section, we give a uniform law of large numbers for ergodic processes, under a smoothness assumption. The proof is standard, so it is omitted. [See, e.g., Theorem 2.4.1 of van der Vaart and Wellner (1996) for the idea, or see Nishiyama (2009).]

THEOREM 3.1. *Let $(E, \mathcal{E})$ be a measurable space. Let $\Theta$ be a set which is totally bounded with respect to the semimetric $\rho$. Let a family $\{f(\cdot; \theta); \theta \in \Theta\}$ of measurable functions on $E$ be given. Suppose that there exists a measurable function $K$ such that*

(3) $$|f(x; \theta) - f(x; \theta')| \leq K(x)\rho(\theta, \theta') \qquad \forall \theta, \theta' \in \Theta.$$

(i) *Suppose that the $E$-valued random process $\{X_t\}_{t \in [0, \infty)}$ is ergodic with the invariant law $\mu$, that is, for any $\mu$-integrable function $g$*
$$\frac{1}{T} \int_0^T g(X_t)\, dt \xrightarrow{p} \int_E g(x)\mu(dx).$$
*If all $f(\cdot; \theta)$ and $K$ are $\mu$-integrable, then*
$$\sup_{\theta \in \Theta} \left| \frac{1}{T} \int_0^T f(X_t; \theta)\, dt - \int_E f(x; \theta)\mu(dx) \right| = o_{P^*}(1).$$

(ii) *Suppose that the $E$-valued random process $\{X_i\}_{i=1,2,\ldots}$ is ergodic with the invariant law $\mu$, that is, for any $\mu$-integrable function $g$,*
$$\frac{1}{n} \sum_{i=1}^n g(X_i) \xrightarrow{p} \int_E g(x)\mu(dx).$$
*If all $f(\cdot; \theta)$ and $K$ are $\mu$-integrable, then*
$$\sup_{\theta \in \Theta} \left| \frac{1}{n} \sum_{i=1}^n f(X_i; \theta) - \int_E f(x; \theta)\mu(dx) \right| = o_{P^*}(1).$$

REMARK. The smoothness assumption (3) can be replaced by "bracketing." See Theorem 2.4.1 of van der Vaart and Wellner (1996).

## 4. Ergodic diffusion processes.

4.1. *Regularity conditions.* Let us consider the diffusion process model introduced in the first paragraph of Section 1. We shall list up some conditions. We suppose that there exists a parametric family of $d$-dimensional vector-valued functions $\{\dot{S}(\cdot;\theta); \theta \in \Theta\}$ on $I$ which satisfies the following conditions. Typically, they may be considered to be the derivatives of $S(\cdot;\theta)$ with respect to $\theta$, that is, $\dot{S}(\cdot;\theta) = (\frac{\partial}{\partial \theta_1}S(\cdot;\theta),\ldots,\frac{\partial}{\partial \theta_d}S(\cdot;\theta))^T$. The function $\Lambda$ appearing in A1 and A3 may be chosen to be common without loss of generality.

A1. $\Theta$ is a compact subset of $\mathbb{R}^d$. There exists a measurable function $\Lambda$ on $I$ such that at the true $\theta_0 \in \Theta$,
$$S(x;\theta) - S(x;\theta_0) = \dot{S}(x;\theta_0)^T(\theta - \theta_0) + \Lambda(x)\epsilon(x;\theta,\theta_0),$$
where $\sup_{x \in I} |\epsilon(x;\theta,\theta_0)| = o(\|\theta - \theta_0\|)$ as $\theta \to \theta_0$.

A2. There exists a constant $K > 0$ such that
$$\sup_{\theta \in \Theta} |S(x;\theta) - S(x';\theta)| \leq K|x - x'|;$$
$$\sup_{\theta \in \Theta} \|\dot{S}(x;\theta) - \dot{S}(x';\theta)\| \leq K|x - x'|;$$
$$\sup_{h \in H} |\sigma^2(x;h) - \sigma^2(x';h)| \leq K|x - x'|.$$

A3. There exists a measurable function $\Lambda$ on $I$ such that
$$\sup_{\theta \in \Theta} |S(x;\theta)| \leq \Lambda(x);$$
$$\sup_{\theta \in \Theta} \|\dot{S}(x;\theta)\| \leq \Lambda(x);$$
$$c := \inf_{h \in H} \inf_{x \in I} \sigma^2(x;h) > 0;$$
$$\|\dot{S}(x;\theta) - \dot{S}(x;\theta)\| \leq \Lambda(x)\|\theta - \theta'\| \qquad \forall \theta, \theta' \in \Theta;$$
$$|\sigma^2(x;h) - \sigma^2(x;h')| \leq \Lambda(x)\,d_H(h,h') \qquad \forall h, h' \in H.$$

A4. $\sup_{t \in \mathbb{R}} E(\Lambda(X_t)^8 + |X_t|^4) < \infty$.

A5. The process $X = (X_t)_{t \in [0,\infty)}$ is ergodic. We denote by $\mu$ the invariant measure under the true $(\theta_0, h_0)$, and we assume that it satisfies $\int_I \Lambda(x)^2(1+|x|)\mu(dx) < \infty$.

A6. The matrix
$$I(\theta_0, h_0) = \int_I \frac{\dot{S}(x;\theta_0)\dot{S}(x;\theta_0)^T}{\sigma^2(x;h_0)}\mu(dx)$$
is invertible.

A7. The metric entropy condition for $(H, d_H)$ is satisfied:

$$\int_0^1 \sqrt{\log N(H, d_H, \varepsilon)}\, d\varepsilon < \infty.$$

A8. For every $\varepsilon > 0$,

$$\inf_{\theta\, :\, \|\theta - \theta_0\| > \varepsilon} \left\| \int_I \frac{\dot{S}(x;\theta)}{\sigma^2(x;h_0)} [S(x;\theta_0) - S(x;\theta)] \mu(dx) \right\| > 0.$$

REMARK. The last assumption in A2 implies that $\sigma^2(x;h_0) \le C(1 + |x|)$ for a constant $C > 0$.

To close this subsection, let us discuss the possibility of the choice of the nuisance parameter space $(H, d_H)$.

EXAMPLE 1 (Parametric model). When $(H, d_H)$ is a compact subset of a finite-dimensional Euclidean space, the metric entropy condition A7 is indeed satisfied. So the main restriction is the Lipschitz continuity of $h \mapsto \sigma^2(\cdot; h)$ in A3. This situation is more general than that in Yoshida (1992) and Kessler (1997), although, as announced in Section 1, our result does not include theirs.

EXAMPLE 2 (The class of smooth functions). Let us consider the parametrization $\sigma(x;h) = h(x)$ where $h$ is an element of the class $H = C_M^\alpha(I)$ defined below. We equip the function space $H$ with the uniform metric $\|\cdot\|_\infty$ for which the last requirement in A3 is always fulfilled. To check A7, first we consider the case where $I$ is a bounded subset of $\mathbb{R}$, and next we give some remarks for the general case.

We take the material below from Section 2.7.1 of van der Vaart and Wellner (1996). Let $I$ be a bounded, convex subset of $\mathbb{R}^q$. (In the current example of one-dimensional diffusions, we are considering the case $q = 1$, but for the generality we set $q$ to be a general positive integer; see Section 5.) Let $\alpha > 0$ and $M > 0$ be given, and let $\underline{\alpha}$ be the greatest integer smaller than $\alpha$. For any vector $k = (k_1, \ldots, k_q)$ of $q$ integers, we define

$$D^k = \frac{\partial^{k\cdot}}{\partial x_1^{k_1} \cdots \partial x_q^{k_q}},$$

where $k_\cdot = \sum_{i=1}^q k_i$. We denote by $C_M^\alpha(I)$ the class of functions defined on $I$ such that

$$\max_{k_\cdot \le \underline{\alpha}} \sup_{\mathbf{x}} |D^k h(\mathbf{x})| + \max_{k_\cdot = \underline{\alpha}} \sup_{\mathbf{x},\mathbf{y}} \frac{|D^k h(\mathbf{x}) - D^k h(\mathbf{y})|}{\|\mathbf{x} - \mathbf{y}\|^{\alpha - \underline{\alpha}}} \le M,$$

where the sumprema are taken over all $\mathbf{x}, \mathbf{y}$ in the interior of $I$ with $\mathbf{x} \neq \mathbf{y}$. Then there exists a constant $K > 0$ depending only on $\alpha$ and $q$, such that

$$\log N(C_M^\alpha(I), \|\cdot\|_\infty, \varepsilon) \leq K\lambda(I^1)\left(\frac{M}{\varepsilon}\right)^{q/\alpha},$$

where $\lambda(I^1)$ is the Lebesgue measure of the set $\{\mathbf{x} : \|\mathbf{x} - I\| < 1\}$. Hence, the metric entropy condition A7 is satisfied if $q/(2\alpha) < 1$, and therefore our theory works.

When $I = \mathbb{R}^q$, we shall restrict out attention, for example, to the following class $H$ of functions on $\mathbb{R}^q$. Let $I_0$ be a bounded, convex subset of $\mathbb{R}^q$, and we suppose that the restriction of $h \in H$ to $I_0$ belongs to $C_M(I_0)$ and that

$$(4) \qquad \sup_{\mathbf{x}\in\mathbb{R}^q} |h(\mathbf{x}) - h'(\mathbf{x})| \leq L \sup_{\mathbf{x}\in I_0} |h(\mathbf{x}) - h'(\mathbf{x})| \qquad \forall h, h' \in H,$$

for a constant $L > 0$. Then both the last condition of A3 and A7 are satisfied. The condition (4) is satisfied if we assume, for example, either of the following:

(i) $h$ is known on $I_0^c$;
(ii) when $q = 1$ and $I_0 = [l_0, r_0]$, each $h$ is constant on $(-\infty, l_0]$ and on $[r_0, \infty)$.

Although the examples (i) and (ii) might look restrictive, it should be noted that in practice we can choose an arbitrary large $I_0$.

Another way to deal with the unbounded case $I = \mathbb{R}^q$ is to consider the parametrization

$$\sigma(\mathbf{x}; h) = h(u(\mathbf{x})), \qquad h \in C_M^\alpha(I_0),$$

where $I_0$ being a bounded, convex subset of $\mathbb{R}^{q'}$ and $u : \mathbb{R}^q \to I_0$ is a fixed function. If $q'/(2\alpha) < 1$, then both the last requirement in A3 and the metric entropy condition A7 are satisfied for the uniform metric $d_H = \|\cdot\|_\infty$.

Instead of $C_M^\alpha(I)$, another possibility of the choice of $H$ which satisfies the metric entropy condition for the uniform metric is the Sobolev class; see Example 19.10 of van der Vaart (1998).

4.2. *Results.* As announced in Section 1, we propose to use the estimating function

$$\Psi_n(\theta, h) = \frac{1}{t_n^n} \sum_{i=1}^n \frac{\dot{S}(X_{t_{i-1}^n}; \theta)}{\sigma^2(X_{t_{i-1}^n}; h)} [X_{t_i^n} - X_{t_{i-1}^n} - S(X_{t_{i-1}^n}; \theta)|t_i^n - t_{i-1}^n|],$$

whose compensator is

$$\widetilde{\Psi}_n(\theta, h) = \frac{1}{t_n^n} \sum_{i=1}^n \frac{\dot{S}(X_{t_{i-1}^n}; \theta)}{\sigma^2(X_{t_{i-1}^n}; h)} \left[\int_{t_{i-1}^n}^{t_i^n} S(X_t, \theta_0)\, dt - S(X_{t_{i-1}^n}; \theta)|t_i^n - t_{i-1}^n|\right].$$

Then we have the following two lemmas.

LEMMA 4.1.  *Assume $\Delta_n \to 0$ and $t_n^n \to \infty$. Equip the space $\Theta \times H$ with the metric $\rho = \|\cdot\| \vee d_H$. Under A1–A5 and A7, $\sqrt{t_n^n}(\Psi_n - \widetilde{\Psi}_n)$ converges weakly in $C_\rho(\Theta \times H)$ to a zero-mean Gaussian process $Z$ with the covariance*

$$EZ(\theta,h)Z(\theta',h')^T = \int_I \frac{\dot{S}(x;\theta)\dot{S}(x;\theta')^T}{\sigma^2(x;h)\sigma^2(x;h')}\sigma^2(x;h_0)\mu(dx).$$

*In particular, the random variable $Z(\theta_0, h_0)$ is distributed with $\mathcal{N}(0, I(\theta_0, h_0))$.*

LEMMA 4.2.  *Assume $\Delta_n = o((t_n^n)^{-1})$ and $t_n^n \to \infty$. Under A1–A6, for any random sequence $(\theta_n, h_n)$ such that $\|\theta_n - \theta_0\| \vee d_H(h_n, h_0) = o_{P^*}(1)$, it holds that*

$$\widetilde{\Psi}_n(\theta_n, h_n) - \widetilde{\Psi}_n(\theta_0, h_0) - (-I(\theta_0, h_0))(\theta_n - \theta_0) = o_{P^*}((t_n^n)^{-1/2} + \|\theta_n - \theta_0\|).$$

Combining these lemmas with Theorem 2.1, and noting also $\widetilde{\Psi}_n(\theta_0, h_0) = O_P(\Delta_n^{1/2})$ which will be proved by using Lemma 4.5 below, we can conclude the following theorem.

THEOREM 4.3.  *Assume $\Delta_n = o((t_n^n)^{-1})$ and $t_n^n \to \infty$. Under A1–A7, for any random sequence $(\widehat{\theta}_n, \widehat{h}_n)$, such that*

$$\|\widehat{\theta}_n - \theta_0\| = o_{P^*}(1), \qquad d_H(\widehat{h}_n, h_0) = o_{P^*}(1)$$

*and*

$$\Psi_n(\widehat{\theta}_n, \widehat{h}_n) = o_{P^*}((t_n^n)^{-1/2}),$$

*the estimator $\widehat{\theta}_n$ is asymptotically normal and efficient:*

$$\sqrt{t_n^n}(\widehat{\theta}_n - \theta_0) \xrightarrow{d} \mathcal{N}(0, I(\theta_0, h_0)^{-1}).$$

*When A8 is also satisfied, the assumption "$\|\widehat{\theta}_n - \theta_0\| = o_{P^*}(1)$" is automatically satisfied.*

In the above theorem, the only assumption which we cannot check in the course of computing the data is the consistency "$d_H(\widehat{h}_n, h_0) = o_{P^*}(1)$," because it involves the true value $h_0$ of the unknown parameter $h \in H$. When $\{\sigma^2(\cdot; h); h \in H\}$ is a class of functions $\sigma^2(\cdot; h) = h(\cdot)$ where $H$ is a class of smooth functions, one may think that a kernel estimator is a candidate for $\widehat{h}_n$. As stated above, in view of the Lipschitz condition of $h \mapsto \sigma^2(\cdot; h)$ (the last condition in A3), it is convenient to consider the consistency with respect to the uniform metric. However, to show the consistency of the kernel estimator with respect to the uniform metric is a task. Generally speaking, showing the consistency for infinite-dimensional parameter is not a trivial problem, which should be solved by independent articles. See, for example, Hoffmann (2001). Below, we give a general way to show the consistency of a least square estimator.

THEOREM 4.4. *Assume $\Delta_n \to 0$ and $t_n^n \to \infty$. Assume A2–A5 and A7. Suppose that*

$$\inf_{h:\, d_H(h,h_0)>\varepsilon} \int_I |\sigma^2(x;h) - \sigma^2(x;h_0)|^2 \mu(dx) > 0 \qquad \forall \varepsilon > 0$$

*is satisfied. If the random element $\widehat{h}_n$ satisfies $\mathcal{A}_n(\widehat{h}_n) \leq \inf_{h \in H} \mathcal{A}_n(h) + o_{P^*}(1)$ where*

$$\mathcal{A}_n(h) = \frac{1}{t_n^n} \sum_{i=1}^{n} \left| \frac{|X_{t_i^n} - X_{t_{i-1}^n}|^2}{|t_i^n - t_{i-1}^n|} - \sigma^2(X_{t_{i-1}^n}, h) \right|^2 |t_i^n - t_{i-1}^n|,$$

*then it holds that $d_H(\widehat{h}_n, h_0) = o_{P^*}(1)$.*

4.3. *Proofs.* Before the proofs, we state a lemma which is well known.

LEMMA 4.5. *Let $X$ be a solution to the SDE (1) for $(\theta, h) = (\theta_0, h_0)$. Assume $|t_i^n - t_{i-1}^n| \leq 1$.*

(i) *For any $k \geq 2$, there exists a constant $C_k > 0$, depending only on $k$, such that*

$$E \sup_{t \in [t_{i-1}^n, t_i^n]} |X_t - X_{t_{i-1}^n}|^k \leq C_k \sup_{s \in \mathbb{R}} E\{|S(X_s; \theta_0)|^k + |\sigma(X_s; h_0)|^k\} |t_i^n - t_{i-1}^n|^{k/2}$$

$$=: D_k |t_i^n - t_{i-1}^n|^{k/2},$$

*provided the right-hand side is finite.*

(ii) *For any $k \geq 2$ and any measurable function $f$, $g$, it holds that*

$$\sup_{t \in [t_{i-1}^n, t_i^n]} E(|X_t - X_{t_{i-1}^n}|^{k/2} |f(X_{t_{i-1}^n})||g(X_t)|)$$

$$\leq (D_k |t_i^n - t_{i-1}^n|^{k/2})^{1/2} \sup_{s \in \mathbb{R}} (E|f(X_s)|^4)^{1/4} \sup_{s \in \mathbb{R}} (E|g(X_s)|^4)^{1/4},$$

*provided the right-hand side is finite.*

PROOF. The assertion (i) is well known. (Use Hölder's inequality and Burkholder–Davis–Gundy's inequality for $\int_{t_{i-1}^n}^{t_i^n} |S(X_s; \theta_0)| \, ds$ and $\sup_{t \in [t_{i-1}^n, t_i^n]} |\int_{t_{i-1}^n}^{t} \sigma(X_s; h_0) \, dW_s|$.) The assertion (ii) follows from Hölder's inequality and (i). □

During the proofs, we write

$$\psi(x; \theta, h) = \frac{\dot{S}(x; \theta)}{\sigma^2(x; h)},$$

which is a $d$-dimensional vector-valued function. For each component $\psi^{(j)}(x;\theta, h)$ $(j=1,\ldots,d)$, it holds that

$$
\begin{aligned}
|\psi^{(j)}&(x;\theta,h) - \psi^{(j)}(x';\theta,h)| \\
&\leq \frac{|\dot{S}^{(j)}(x;\theta) - \dot{S}^{(j)}(x';\theta)|}{\sigma^2(x;h)} \\
&\quad + |\dot{S}^{(j)}(x';\theta)|\left|\frac{1}{\sigma^2(x;h)} - \frac{1}{\sigma^2(x';h)}\right| \\
&\leq \left\{\frac{1}{c} + \frac{\Lambda(x)}{c^2}\right\}K|x - x'|
\end{aligned}
\tag{5}
$$

and that

$$
\begin{aligned}
|\psi^{(j)}&(x;\theta,h) - \psi^{(j)}(x;\theta',h')| \\
&\leq \frac{|\dot{S}^{(j)}(x;\theta) - \dot{S}^{(j)}(x;\theta')|}{\sigma^2(x;h)} + |\dot{S}^{(j)}(x;\theta')|\left|\frac{1}{\sigma^2(x;h)} - \frac{1}{\sigma^2(x;h')}\right| \\
&\leq \frac{|\dot{S}^{(j)}(x;\theta) - \dot{S}^{(j)}(x;\theta')|}{c} + |\dot{S}^{(j)}(x;\theta')|\frac{|\sigma^2(x;h) - \sigma^2(x;h')|}{c^2} \\
&\leq \left\{\frac{\Lambda(x)}{c} + \frac{|\Lambda(x)|^2}{c^2}\right\}(\|\theta - \theta'\| \vee d_H(h,h')).
\end{aligned}
\tag{6}
$$

PROOF OF LEMMA 4.1. We apply Theorem 3.4.2 of Nishiyama (2000b) [or, see Theorem 3.3 of van der Vaart and van Zanten (2005)] to the terminals $M^{n,\theta,h}_{t^n_n}$ of the continuous martingales $t \rightsquigarrow M^{n,\theta,h}_t$ given by

$$M^{n,\theta,h}_t = \frac{1}{\sqrt{t^n_n}}\sum_{i=1}^n \psi(X_{t^n_{i-1}};\theta,h)\int_{t^n_{i-1}\wedge t}^{t^n_i \wedge t} \sigma(X_s;h_0)\,dW_s.$$

For the finite-dimensional convergence, it is sufficient to show the convergence of predictable covariation. This is done as follows.

$$
\begin{aligned}
\langle M^{n,\theta,h}&, M^{n,\theta',h'}\rangle_{t^n_n} \\
&= \frac{1}{t^n_n}\sum_{i=1}^n \psi(X_{t^n_{i-1}};\theta,h)\psi(X_{t^n_{i-1}};\theta',h')^T \int_{t^n_{i-1}}^{t^n_i} \sigma^2(X_s;h_0)\,ds \\
&= \frac{1}{t^n_n}\int_0^{t^n_n} \psi(X_t;\theta,h)\psi(X_t;\theta',h')^T \sigma^2(X_t;h_0)\,dt + o_P(1) \\
&\xrightarrow{p} \int_I \psi(x;\theta,h)\psi(x;\theta',h')^T \sigma^2(x;h_0)\mu(dx) \\
&= \int_I \frac{\dot{S}(x;\theta)\dot{S}(x;\theta')^T}{\sigma^2(x;h)\sigma^2(x;h')}\sigma^2(x;h_0)\mu(dx)
\end{aligned}
\tag{7}
$$

$$= I(\theta_0, h_0) \qquad \text{if } (\theta, h) = (\theta', h') = (\theta_0, h_0).$$

Here, to show (7), we have used the bound (5) and Lemma 4.5(ii) twice.

To establish Nishiyama's condition [ME], let us observe the following fact to check the metric entropy condition for the product space $\Theta \times H$.

In general, if $(D, d)$ and $(E, e)$ are two semimetric spaces, then the covering number of the product space $D \times E$ with respect to the maximum semimetric $d \vee e$, namely $N(D \times E, d \vee e, \varepsilon)$, is bounded by $N(D, d, \varepsilon) \cdot N(E, e, \varepsilon)$. To see this claim, let $B_i$, $i = 1, \ldots, N(D, d, \varepsilon)$ be an $\varepsilon$-covering of $D$, and let $C_j$, $j = 1, \ldots, N(E, e, \varepsilon)$ be an $\varepsilon$-covering of $E$. Then the diameters of the sets $B_i \times C_j \subset D \times E$ with respect to $d \vee e$ are smaller than $\varepsilon$, thus these sets form an $\varepsilon$-covering of $D \times E$. The claim has been proved. Consequently, the metric entropy condition

$$\int_0^1 \sqrt{\log N(D \times E, d \vee e, \varepsilon)} \, d\varepsilon < \infty$$

is satisfied if

$$\int_0^1 \sqrt{\log N(D, d, \varepsilon)} \, d\varepsilon < \infty \quad \text{and} \quad \int_0^1 \sqrt{\log N(E, e, \varepsilon)} \, d\varepsilon < \infty.$$

Now, since $\Theta$ is compact with respect to the Euclidean metric, the metric entropy condition for $\Theta$ is satisfied. So, with A7 in hands, the metric entropy condition for the product space $\Theta \times H$ is satisfied, and it remains only to show that the *quadratic modulus* is bounded in probability; that is, the claim that each component of the matrix

$$\sup_{(\theta, h) \neq (\theta', h')} \frac{\langle M^{n,\theta,h} - M^{n,\theta',h'} \rangle_{t_n^n}}{(\|\theta - \theta'\| \vee d(h, h'))^2}$$

is bounded in probability. In view of (6), the absolute value of each component of this matrix is bounded by

$$\frac{1}{t_n^n} \sum_{i=1}^n \int_{t_{i-1}^n}^{t_i^n} \left| \frac{\Lambda(X_{t_{i-1}^n})}{c} + \frac{|\Lambda(X_{t_{i-1}^n})|^2}{c^2} \right|^2 \sigma^2(X_s; h_0) \, ds.$$

The expectation of this random valuable is bounded by

$$\sup_i \sqrt{E \left| \frac{\Lambda(X_{t_{i-1}^n})}{c} + \frac{|\Lambda(X_{t_{i-1}^n})|^2}{c^2} \right|^4} \cdot \sup_s \sqrt{E \sigma^4(X_s; h_0)},$$

which is $O(1)$ by A4. Thus, the quadratic modulus is bounded in probability. □

PROOF OF LEMMA 4.2. It follows from Lemma 4.5 that uniformly in $\theta, h$,

$$\widetilde{\Psi}_n(\theta, h) - \widetilde{\Psi}_n(\theta_0, h_0)$$
$$= \frac{1}{t_n^n} \sum_{i=1}^n \psi(X_{t_{i-1}^n}; \theta, h) \int_{t_{i-1}^n}^{t_i^n} [S(X_t; \theta_0) - S(X_t; \theta)] \, dt + O_P(\Delta_n^{1/2})$$
$$= \frac{1}{t_n^n} \int_0^{t_n} \psi(X_t; \theta, h)[S(X_t; \theta_0) - S(X_t; \theta)] \, dt + O_P(\Delta_n^{1/2}).$$

The remainder term of this approximation is actually $o_P((t_n^n)^{-1/2})$. Furthermore, it holds for any (possibly random) sequence $(\theta_n, h_n)$ converging in outer probability to $(\theta_0, h_0)$ that

$$\frac{1}{t_n^n} \int_0^{t_n^n} \psi(X_t; \theta_n, h_n)[S(X_t; \theta_0) - S(X_t; \theta_n)] \, dt$$
$$= \frac{1}{t_n^n} \int_0^{t_n^n} \psi(X_t; \theta_n, h_n) \dot{S}(X_t; \theta_0)^T \, dt (\theta_0 - \theta_n) + o_{P^*}(\|\theta_n - \theta_0\|)$$
(8)
$$= \frac{1}{t_n^n} \int_0^{t_n^n} \psi(X_t; \theta_0, h_0) \dot{S}(X_t; \theta_0)^T \, dt (\theta_0 - \theta_n) + o_{P^*}(\|\theta_n - \theta_0\|)$$
$$= \int_I \psi(x; \theta_0, h_0) \dot{S}(x; \theta_0)^T \mu(dx) (\theta_0 - \theta_n) + o_{P^*}(\|\theta_n - \theta_0\|)$$
$$= -I(\theta_0, h_0)(\theta_n - \theta_0) + o_{P^*}(\|\theta_n - \theta_0\|).$$

To prove (8) in the above computation, use (6) to show that for every $j, k = 1, \ldots, d$

$$\left| \frac{1}{t_n^n} \int_0^{t_n^n} \psi^{(j)}(X_t; \theta_n, h_n) \dot{S}^{(k)}(X_t; \theta_0) \, dt - \frac{1}{t_n^n} \int_0^{t_n^n} \psi^{(j)}(X_t; \theta_0, h_0) \dot{S}^{(k)}(X_t; \theta_0) \, dt \right|$$
$$\leq \frac{1}{t_n^n} \int_0^{t_n^n} \left\{ \frac{\Lambda(X_t)}{c} + \frac{|\Lambda(X_t)|^2}{c^2} \right\} |\dot{S}^{(k)}(X_t; \theta_0)| \, dt \cdot \|\theta_n - \theta_0\| \vee d_H(h_n, h_0)$$
$$= O_P(1) \cdot o_{P^*}(1)$$
$$= o_{P^*}(1).$$

The proof is complete. □

PROOF OF THEOREM 4.3. By Lemma 4.5, it is easy to see that $\widetilde{\Psi}_n(\theta_0, h_0) = O_P(\Delta_n^{1/2}) = o_P((t_n^n)^{-1/2})$. So the main assertion follows from Theorem 2.1 with help from Lemmas 4.1 and 4.2. On the other hand, since it follows from Lemma 4.5(ii) and Theorem 3.1(i) that $\sup_{\theta, h} \|\Psi_n(\theta, h) - \Psi(\theta, h)\| = o_{P^*}(1)$, where

$$\Psi(\theta, h) = \int_I \frac{\dot{S}(x; \theta)}{\sigma^2(x; h)} [S(x; \theta_0) - S(x; \theta)] \mu(dx),$$

the assertion that the consistency "$\|\widehat{\theta}_n - \theta_0\| = o_{P^*}(1)$" automatically follows from A8 is immediate from Theorem 2.3. □

PROOF OF THEOREM 4.4. Put
$$M_n(h) = \frac{1}{t_n^n} \sum_{i=1}^n |\sigma^2(X_{t_{i-1}^n}; h) - \sigma^2(X_{t_{i-1}^n}; h_0)|^2 |t_i^n - t_{i-1}^n|,$$

$$M(h) = \int_I |\sigma^2(x; h) - \sigma^2(x; h_0)|^2 \mu(dx).$$

Let us apply Corollary 3.2.3 of van der Vaart and Wellner (1996) to the above $M_n$ and $M$ for the given $\widehat{h}_n$ which is the solution of $\mathcal{A}_n(\widehat{h}_n) \leq \inf_{h \in H} \mathcal{A}_n(h) + o_{P^*}(1)$. By Lemma 4.5(ii) and Theorem 3.1(i), it is not difficult to see that $\sup_{h \in H} |M_n(h) - M(h)| = o_{P^*}(1)$, so it is sufficient to show that $M_n(\widehat{h}_n) = o_{P^*}(1)$.

Observe that
$$\frac{1}{t_n^n} \sum_{i=1}^n \left| \frac{|X_{t_i^n} - X_{t_{i-1}^n}|^2}{|t_i^n - t_{i-1}^n|} - \sigma^2(X_{t_{i-1}^n}; h) \right|^2 |t_i^n - t_{i-1}^n|$$
$$= \frac{1}{t_n^n} \sum_{i=1}^n \left| \frac{|X_{t_i^n} - X_{t_{i-1}^n}|^2}{|t_i^n - t_{i-1}^n|} - \sigma^2(X_{t_{i-1}^n}; h_0) \right|^2 |t_i^n - t_{i-1}^n|$$
$$+ \frac{2}{t_n^n} \sum_{i=1}^n \left( \frac{|X_{t_i^n} - X_{t_{i-1}^n}|^2}{|t_i^n - t_{i-1}^n|} - \sigma^2(X_{t_{i-1}^n}; h_0) \right)$$
$$\times (\sigma^2(X_{t_{i-1}^n}; h_0) - \sigma^2(X_{t_{i-1}^n}; h)) |t_i^n - t_{i-1}^n|$$
$$+ \frac{1}{t_n^n} \sum_{i=1}^n |\sigma^2(X_{t_{i-1}^n}; h_0) - \sigma^2(X_{t_{i-1}^n}; h)|^2 |t_i^n - t_{i-1}^n|.$$

Let us prove that the supremum with respect to $h$ of the absolute value of the second term on the right-hand side converges in outer probability to zero [say, the claim (a)].

Since we have from Itô's formula that
$$|X_{t_i^n} - X_{t_{i-1}^n}|^2 = 2 \int_{t_{i-1}^n}^{t_i^n} (X_s - X_{t_{i-1}^n}) S(X_s; \theta_0) \, ds$$
$$+ 2 \int_{t_{i-1}^n}^{t_i^n} (X_s - X_{t_{i-1}^n}) \sigma(X_s; h_0) \, dW_s + \int_{t_{i-1}^n}^{t_i^n} \sigma^2(X_s; h_0) \, ds,$$

it is sufficient to show that $C_{1,n} = o_P(1)$, $\sup_{h \in H} |C_{2,n}(h)| = o_{P^*}(1)$ and $C_{3,n} = o_P(1)$, where
$$C_{1,n} = \frac{1}{t_n^n} \sum_{i=1}^n \left| \int_{t_{i-1}^n}^{t_i^n} (X_s - X_{t_{i-1}^n}) S(X_s; \theta_0) \, ds \right| \Lambda(X_{t_{i-1}^n}),$$

$$C_{2,n}(h) = \frac{1}{t_n^n} \sum_{i=1}^{n} \int_{t_{i-1}^n}^{t_i^n} (X_s - X_{t_{i-1}^n})\sigma(X_s; h_0)\, dW_s(\sigma^2(X_{t_{i-1}^n}; h_0) - \sigma^2(X_{t_{i-1}^n}; h)),$$

$$C_{3,n} = \frac{1}{t_n^n} \sum_{i=1}^{n} \int_{t_{i-1}^n}^{t_i^n} |\sigma^2(X_s; h_0) - \sigma^2(X_{t_{i-1}^n}; h_0)|\, ds \Lambda(X_{t_{i-1}^n}).$$

By using Lemma 4.5(ii), we easily have $EC_{1,n} \to 0$ and $EC_{3,n} \to 0$. On the other hand, by using Theorem 3.4.2 of Nishiyama (2000b), it holds that $C_{2,n}$ converges weakly to zero in $C_{d_H}(H)$ (recall the argument in the proof of Lemma 4.1). Therefore, we have $\sup_{h \in H} |C_{2,n}(h)| = o_{P^*}(1)$.

Hence, the claim (a) is true, and we have that $\mathcal{A}_n(\widehat{h}_n) \leq \inf_{h \in H} \mathcal{A}_n(h) + o_{P^*}(1)$ implies that $M_n(\widehat{h}_n) = o_{P^*}(1)$. The proof is finished. □

## 5. Ergodic time series.

5.1. *Model and regularity conditions.* Let us consider the time series model given by

$$X_i = \widetilde{S}(X_{i-1}, \ldots, X_{i-q_1}; \theta) + \widetilde{\sigma}(X_{i-1}, \ldots, X_{i-q_2}; h)w_i.$$

By putting $q = q_1 \vee q_2$ and changing the domains of the functions $\widetilde{S}$ and $\widetilde{\sigma}$, without loss of generality, we can write

$$X_i = S(\mathbf{X}_i; \theta) + \sigma(\mathbf{X}_i; h)w_i,$$

where $\mathbf{X}_i = (X_{i-1}, \ldots, X_{i-q})$ and $S(\cdot; \theta)$ and $\sigma(\cdot; h)$ are some measurable functions on $\mathbb{R}^q$. For simplicity, we assume that the initial values $(X_0, \ldots, X_{1-q}) = (x_0, \ldots, x_{1-q})$ are fixed.

As for the noise $\{w_i\}$, we consider the following two cases:

CASE G (Gaussian). $\{w_i\}$ are independently, identically distributed with $\mathcal{N}(0, 1)$.

CASE M (Martingale). $E[w_i | \mathcal{F}_{i-1}] = 0$ and $E[w_i^2 | \mathcal{F}_{i-1}] = 1$ almost surely, where $\mathcal{F}_i = \sigma\{X_j : j \leq i\}$.

Clearly, the Case G is a special case of the Case M. When we do not especially declare the restriction to the Case G, we consider the Case M in principle.

Let us list up some conditions which have the same fashion as those in Section 4.1. We suppose that there exists a parametric family of $d$-dimensional vector-valued functions $\{\dot{S}(\cdot; \theta); \theta \in \Theta\}$ on $\mathbb{R}^q$ which satisfies the following conditions. Typically, they may be considered to be the derivatives of $S(\cdot; \theta)$ with respect to $\theta$, that is, $\dot{S}(\cdot; \theta) = (\frac{\partial}{\partial \theta_1} S(\cdot; \theta), \ldots, \frac{\partial}{\partial \theta_d} S(\cdot; \theta))^T$. The function $\Lambda$ appearing in B1 and B2 may be chosen to be common without loss of generality.

B1. $\Theta$ is a compact subset of $\mathbb{R}^d$. There exists a measurable function $\Lambda$ on $\mathbb{R}^q$ such that at the true $\theta_0 \in \Theta$,

$$S(\mathbf{x};\theta) - S(\mathbf{x};\theta_0) = \dot{S}(\mathbf{x};\theta_0)^T(\theta - \theta_0) + \Lambda(\mathbf{x})\epsilon(\mathbf{x};\theta,\theta_0),$$

where $\sup_{\mathbf{x} \in \mathbb{R}^q} |\epsilon(\mathbf{x};\theta,\theta_0)| = o(\|\theta - \theta_0\|)$ as $\theta \to \theta_0$.

B2. There exists a measurable function $\Lambda$ on $\mathbb{R}^q$ such that

$$\sup_{\theta \in \Theta} |S(\mathbf{x};\theta)| \leq \Lambda(\mathbf{x});$$

$$\sup_{\theta \in \Theta} \|\dot{S}(\mathbf{x};\theta)\| \leq \Lambda(\mathbf{x});$$

$$\sigma^2(\mathbf{x};h_0) \leq \Lambda(\mathbf{x}), \qquad c := \inf_{h \in H} \inf_{\mathbf{x} \in \mathbb{R}^q} \sigma^2(\mathbf{x};h) > 0;$$

$$\|\dot{S}(\mathbf{x};\theta) - \dot{S}(\mathbf{x};\theta')\| \leq \Lambda(\mathbf{x})\|\theta - \theta'\| \qquad \forall \theta, \theta' \in \Theta;$$

$$|\sigma^2(\mathbf{x};h) - \sigma^2(\mathbf{x};h')| \leq \Lambda(\mathbf{x})\, d_H(h,h') \qquad \forall h, h' \in H.$$

B3. The process $\{X_i\}_{i=1,2,\ldots}$ is ergodic under the true $(\theta_0, h_0)$ in the sense that for $q' = q$ and $q+1$ there exists the invariant measure $\mu_{q'}$ such that for every $\mu_{q'}$-integrable function $f$

$$\frac{1}{n}\sum_{i=1}^n f(X_{i-1},\ldots,X_{i-q'}) \xrightarrow{p} \int_{\mathbb{R}^{q'}} f(x_1,\ldots,x_{q'})\mu_{q'}(dx_1\cdots dx_{q'}).$$

We also assume that

$$\int_{\mathbb{R}^q} \Lambda(\mathbf{x})^5 \mu_q(d\mathbf{x}) < \infty,$$

$$\int_{\mathbb{R}^{q+1}} ||x_0| + \Lambda(x_1,\ldots,x_q)|^4 \mu_{q+1}(dx_0\, dx_1\cdots dx_q) < \infty.$$

B4. The matrix

$$I(\theta_0, h_0) = \int_{\mathbb{R}^q} \frac{\dot{S}(\mathbf{x};\theta_0)\dot{S}(\mathbf{x};\theta_0)^T}{\sigma^2(\mathbf{x};h_0)} \mu_q(d\mathbf{x})$$

is invertible.

B5. The metric entropy condition for $(H, d_H)$ is satisfied:

$$\int_0^1 \sqrt{\log N(H, d_H, \varepsilon)}\, d\varepsilon < \infty.$$

B6. For every $\varepsilon > 0$,

$$\inf_{\theta:\|\theta - \theta_0\| > \varepsilon} \left\| \int_{\mathbb{R}^q} \frac{\dot{S}(\mathbf{x};\theta)}{\sigma^2(\mathbf{x};h_0)}[S(\mathbf{x};\theta_0) - S(\mathbf{x};\theta)]\mu_q(d\mathbf{x}) \right\| > 0.$$

See the end of Section 4.1 for the discussion of the choice of $(H, d_H)$.

5.2. *Results.* In order to explain the idea of our estimating function, let us first consider the Case G. We denote by $P_{n,u}$ the distribution of $\{X_1, \ldots, X_n\}$ under $\theta = \theta_0 + n^{-1/2}u$ and $h = h_0$, where $u \in \mathbb{R}^d$. By an easy computation, the log-likelihood ratio is given by

$$\log \frac{dP_{n,u}}{dP_{n,0}}(X_1, \ldots, X_n)$$

(9)
$$= -\sum_{i=1}^{n} \frac{1}{2\sigma^2(\mathbf{X}_i; h_0)}$$
$$\times \{|X_i - S(\mathbf{X}_i; \theta_0 + n^{-1/2}u)|^2 - |X_i - S(\mathbf{X}_i; \theta_0)|^2\}$$
$$= \Delta_{n,u} - B_{n,u},$$

where

$$\Delta_{n,u} = \sum_{i=1}^{n} \frac{1}{\sigma^2(\mathbf{X}_i; h_0)} \{X_i - S(\mathbf{X}_i; \theta_0)\} \{S(\mathbf{X}_i; \theta_0 + n^{-1/2}u) - S(\mathbf{X}_i; \theta_0)\}$$

and

$$B_{n,u} = \sum_{i=1}^{n} \frac{1}{2\sigma^2(\mathbf{X}_i; h_0)} \{S(\mathbf{X}_i; \theta_0 + n^{-1/2}u) - S(\mathbf{X}_i; \theta_0)\}^2.$$

Under the above regularity conditions, we have

$$\Delta_{n,u} \xrightarrow{d} \mathcal{N}(0, u^T I(\theta_0, h_0) u) \quad \text{and} \quad B_{n,u} \xrightarrow{p} \tfrac{1}{2} u^T I(\theta_0, h_0) u.$$

So it follows from the theory of the *local asymptotic normality* that the distribution of the asymptotically efficient bound in the Case G when $h_0$ is *known* is $\mathcal{N}(0, I(\theta_0, h_0)^{-1})$. That is, if we obtain an estimator $\tilde{\theta}_n$ such that $\sqrt{n}(\tilde{\theta}_n - \theta_0) \xrightarrow{d} \mathcal{N}(0, I(\theta_0, h_0)^{-1})$, it is asymptotically efficient in the sense of the local asymptotic minimax theorem. [See, e.g., Chapter 3.11 of van der Vaart and Wellner (1996).] If the parameter $h$ is *unknown*, then the estimation problem for $\theta$ becomes more difficult. So if we have an estimator which asymptotically behaves as stated above, then we may say that it is asymptotically efficient with the nuisance parameter $h$. This argument is not true in the Case M where the log-likelihood does not equal the formula (9), but we propose to use it for deriving an estimating equation which yields the same asymptotic distribution as the Case G.

Not only in the Case G but also in the case M, differentiating (9) formally, and replacing the true $h_0$ by the unknown parameter $h$, we propose the estimating function

$$\Psi_n(\theta, h) = \frac{1}{n} \sum_{i=1}^{n} \frac{\dot{S}(\mathbf{X}_i; \theta)}{\sigma^2(\mathbf{X}_i; h)} (X_i - S(\mathbf{X}_i; \theta)).$$

Its compensator is

$$\widetilde{\Psi}_n(\theta, h) = \frac{1}{n}\sum_{i=1}^n \frac{\dot{S}(\mathbf{X}_i;\theta)}{\sigma^2(\mathbf{X}_i;h)}(S(\mathbf{X}_i;\theta_0) - S(\mathbf{X}_i;\theta)).$$

Thus, it holds that

$$\sqrt{n}(\Psi_n(\theta,h) - \widetilde{\Psi}_n(\theta,h)) = \frac{1}{\sqrt{n}}\sum_{i=1}^n \frac{\dot{S}(\mathbf{X}_i;\theta)}{\sigma^2(\mathbf{X}_i;h)}\sigma(\mathbf{X}_i;h_0)w_i,$$

which is the summation of a $C_\rho(\Theta \times H)$-valued martingale difference array where $\rho = \|\cdot\| \vee d_H$. By using Jain–Marcus' central limit theorem for martingales given by Nishiyama (1996, 2000a, 2000b), we have the following lemma which plays a key role in our approach.

LEMMA 5.1. *Under* B1–B3 *and* B5, *the sequence of random fields* $\sqrt{n}(\Psi_n - \widetilde{\Psi}_n)$, *with parameter* $(\theta, h)$, *converges weakly in* $C_\rho(\Theta \times H)$ *to a zero-mean Gaussian random field* $Z$ *with the covariance*

$$EZ(\theta,h)Z(\theta',h')^T = \int_{\mathbb{R}^q} \frac{\dot{S}(\mathbf{x};\theta)\dot{S}(\mathbf{x};\theta')^T}{\sigma^2(\mathbf{x};h)\sigma^2(\mathbf{x};h')}\sigma^2(\mathbf{x};h_0)\mu_q(d\mathbf{x}).$$

*In particular, the random variable* $Z(\theta_0, h_0)$ *is distributed with* $\mathcal{N}(0, I(\theta_0, h_0))$.

Another lemma which is necessary to apply Theorem 2.1 is the following.

LEMMA 5.2. *Under* B1–B4, *for any random sequence* $(\theta_n, h_n)$ *such that* $\|\theta_n - \theta_0\| \vee d_H(h_n, h_0) = o_{P^*}(1)$, *it holds that*

$$\widetilde{\Psi}_n(\theta_n, h_n) - \widetilde{\Psi}_n(\theta_0, h_0) - (-I(\theta_0, h_0))(\widehat{\theta}_n - \theta_0) = o_{P^*}(\|\theta_n - \theta_0\|).$$

Noting also that $\tilde{\Psi}_n(\theta_0, h_0) = 0$, we can apply Theorem 2.1 to conclude the following theorem.

THEOREM 5.3. *Under* B1–B5, *for any random sequence* $(\widehat{\theta}_n, \widehat{h}_n)$ *such that*

$$\|\widehat{\theta}_n - \theta_0\| = o_{P^*}(1), \qquad d_H(\widehat{h}_n, h_0) = o_{P^*}(1) \quad and \quad \Psi_n(\widehat{\theta}_n, \widehat{h}_n) = o_{P^*}(n^{-1/2}),$$

*the estimator* $\widehat{\theta}_n$ *is asymptotically normal:*

$$\sqrt{n}(\widehat{\theta}_n - \theta_0) \xrightarrow{d} \mathcal{N}(0, I(\theta_0, h_0)^{-1}).$$

*In particular, in the Case* G, *the estimator* $\widehat{\theta}_n$ *is asymptotically efficient. When* B6 *is also satisfied, the assumption "*$\|\widehat{\theta}_n - \theta_0\| = o_{P^*}(1)$*" is automatically satisfied.*

By the same reason as in Section 4, it is necessary to develop a procedure to construct a consistent estimator $\widehat{h}_n$ for the nuisance parameter $h \in H$. The following theorem gives us an answer.

THEOREM 5.4. *Assume* B2, B3 *and* B5.
(First step: initial estimator for $\theta_0$.) *Suppose the identifiability condition*

$$\inf_{\theta \,:\, \|\theta-\theta_0\|>\varepsilon} \int_{\mathbb{R}^q} |S(\mathbf{x};\theta) - S(\mathbf{x};\theta_0)|^2 \mu_q(d\mathbf{x}) > 0 \qquad \forall \varepsilon > 0$$

*is satisfied. If a random sequence $\theta_n^{LS}$ satisfies $\mathcal{A}_n(\theta_n^{LS}) \leq \inf_{\theta \in \Theta} \mathcal{A}_n(\theta) + o_{P^*}(1)$, where*

$$\mathcal{A}_n(\theta) = \frac{1}{n} \sum_{i=1}^n |X_i - S(\mathbf{X}_i;\theta)|^2,$$

*then it holds that $\|\theta_n^{LS} - \theta_0\| = o_{P^*}(1)$.*

(Second step: consistent estimator for $h_0$.) *Suppose the identifiability condition*

$$\inf_{h \,:\, d_H(h,h_0)>\varepsilon} \int_{\mathbb{R}^q} |\sigma^2(\mathbf{x};h) - \sigma^2(\mathbf{x};h_0)|^2 \mu_q(d\mathbf{x}) > 0 \qquad \forall \varepsilon > 0$$

*is satisfied. Merely by a technical reason, assume that there exists a constant $L_4 > 0$ such that $E[w_i^4|\mathcal{F}_{i-1}] < L_4$ almost surely for all $i$. Using $\theta_n^{LS}$ as above, we define*

$$\mathcal{B}_n(h) = \frac{1}{n} \sum_{i=1}^n ||X_i - S(\mathbf{X}_i;\theta_n^{LS})|^2 - \sigma^2(\mathbf{X}_i;h)|^2.$$

*If a random sequence $\widehat{h}_n$ satisfies $\mathcal{B}_n(\widehat{h}_n) \leq \inf_{h \in H} \mathcal{B}_n(h) + o_{P^*}(1)$, then it holds that $d_H(\widehat{h}_n, h_0) = o_{P^*}(1)$.*

5.3. *Proofs.*

PROOF OF LEMMA 5.1. To show the finite-dimensional convergence is easy. Notice that

$$\begin{aligned}
\left| \frac{\dot{S}(\mathbf{x};\theta)}{\sigma^2(\mathbf{x};h)} - \frac{\dot{S}(\mathbf{x};\theta')}{\sigma^2(\mathbf{x};h')} \right| & \\
&\leq \frac{|\dot{S}(\mathbf{x};\theta) - \dot{S}(\mathbf{x};\theta')|}{\sigma^2(\mathbf{x};h)} + \dot{S}(\mathbf{x}_i;\theta') \left| \frac{1}{\sigma^2(\mathbf{x};h)} - \frac{1}{\sigma^2(\mathbf{x};h')} \right| \\
&\leq \frac{\Lambda(\mathbf{x})}{c} \|\theta - \theta'\| + \frac{\Lambda(\mathbf{x})^2}{c^2} d_H(h,h') \\
&\leq \left\{ \frac{\Lambda(\mathbf{x})}{c} + \frac{\Lambda(\mathbf{x})^2}{c^2} \right\} (\|\theta - \theta'\| \vee d_H(h,h')).
\end{aligned}$$

(10)

The assertion follows from Proposition 4.5 of Nishiyama (2000a). [Or, see Nishiyama (1996) which is easier to read.] □

PROOF OF LEMMA 5.2. Notice that $\widetilde{\Psi}_n(\theta_0, h_0) = 0$. For any (possibly random) sequence $(\theta_n, h_n)$ converging in outer probability to $(\theta_0, h_0)$, it holds that

$$
\begin{aligned}
\widetilde{\Psi}_n(\theta_n, h_n) &- \widetilde{\Psi}_n(\theta_0, h_0) \\
&= \frac{1}{n} \sum_{i=1}^n \frac{\dot{S}(\mathbf{X}_i; \theta_n)}{\sigma^2(\mathbf{X}_i; h_n)} [S(\mathbf{X}_i; \theta_0) - S(\mathbf{X}_i; \theta_n)] \\
&= \frac{1}{n} \sum_{i=1}^n \frac{\dot{S}(\mathbf{X}_i; \theta_n)}{\sigma^2(\mathbf{X}_i; h_n)} \dot{S}(\mathbf{X}_i; \theta_0)(\theta_0 - \theta_n) + o_{P^*}(\|\theta_n - \theta_0\|) \\
&= \frac{1}{n} \sum_{i=1}^n \frac{\dot{S}(\mathbf{X}_i; \theta_0)}{\sigma^2(\mathbf{X}_i; h_0)} \dot{S}(\mathbf{X}_i; \theta_0)(\theta_0 - \theta_n) + o_{P^*}(\|\theta_n - \theta_0\|) \\
&= -I(\theta_0, h_0)(\theta_n - \theta_0) + o_{P^*}(\|\theta_n - \theta_0\|).
\end{aligned}
$$

(11)

To show (11) above, do the same argument as the proof of Lemma 4.2 using (10) instead of (6). □

PROOF OF THEOREM 5.3. The assertions follow from Theorems 2.1 and 2.3 by using also Lemmas 5.1 and 5.2, and Theorem 3.1(ii), respectively. □

PROOF OF THEOREM 5.4. To prove the first step, we will apply Corollary 3.2.3 of van der Vaart and Wellner (1996). We can write $\mathcal{A}_n(\theta) = T_{1,n} + T_{2,n}(\theta) + T_{3,n}(\theta)$ where

$$T_{1,n} = \frac{1}{n} \sum_{i=1}^n |X_i - S(\mathbf{X}_i; \theta_0)|^2,$$

$$T_{2,n}(\theta) = \frac{1}{n} \sum_{i=1}^n (S(\mathbf{X}_i; \theta_0) - S(\mathbf{X}_i; \theta))\sigma(X_i; h_0)w_i,$$

$$T_{3,n}(\theta) = \frac{1}{n} \sum_{i=1}^n |S(\mathbf{X}_i; \theta_0) - S(\mathbf{X}_i; \theta)|^2.$$

The term $T_{1,n}$ converges in probability to a constant $C_1$. On the other hand, by using Proposition 4.5 of Nishiyama (2000a), we have that $\sqrt{n}T_{2,n}$ converges weakly in $C(\Theta)$ to a tight law, thus, $\sup_{\theta \in \Theta} |T_{2,n}(\theta)|$ converges in outer probability to zero. Finally, by Theorem 3.1, it holds that $\sup_{\theta \in \Theta} |T_{3,n}(\theta) - T_3(\theta)| = o_{P^*}(1)$ where

$$T_3(\theta) = \int_{\mathbb{R}^q} |S(\mathbf{x}; \theta_0) - S(\mathbf{x}; \theta)|^2 \mu(d\mathbf{x}).$$

Hence, we have $\sup_{\theta\in\Theta}|\mathcal{A}_n(\theta) - (C_1 + T_3(\theta))| = o_{P^*}(1)$, and van der Vaart and Wellner's (1996) consistency theorem yields the assertion of the first step.

To prove the second step, we shall apply Corollary 3.2.3 of van der Vaart and Wellner (1996) again. Let us first see that $\sup_{h\in H}|\mathcal{B}_n(h) - \widetilde{\mathcal{B}}_n(h)| = o_{P^*}(1)$ where

$$\widetilde{\mathcal{B}}_n(h) = \frac{1}{n}\sum_{i=1}^n ||X_i - S(\mathbf{X}_i;\theta_0)|^2 - \sigma^2(\mathbf{X}_i;h)|^2.$$

Now, notice that

$$\begin{aligned}
&\left| ||X_i - S(\mathbf{X}_i;\theta_n^{LS})|^2 - \sigma^2(\mathbf{X}_i;h)|^2 - ||X_i - S(\mathbf{X}_i;\theta_0)|^2 - \sigma^2(\mathbf{X}_i;h)|^2 \right| \\
&\quad \leq ||X_i - S(\mathbf{X}_i;\theta_n^{LS})|^2 + |X_i - S(\mathbf{X}_i;\theta_0)|^2 - 2\sigma^2(\mathbf{X}_i;h)| \\
&\qquad \times ||X_i - S(\mathbf{X}_i;\theta_n^{LS})|^2 - |X_i - S(\mathbf{X}_i;\theta_0)|^2| \\
&\quad \leq ||X_i - S(\mathbf{X}_i;\theta_n^{LS})|^2 + |X_i - S(\mathbf{X}_i;\theta_0)|^2 - 2\sigma^2(\mathbf{X}_i;h)| \\
&\qquad \times |2X_i - S(\mathbf{X}_i;\theta_n^{LS}) - S(\mathbf{X}_i;\theta_0)| \\
&\qquad \times |S(\mathbf{X}_i;\theta_n^{LS}) - S(\mathbf{X}_i;\theta_0)| \\
&\quad \leq L(X_i, \mathbf{X}_i)\|\theta_n^{LS} - \theta_0\|,
\end{aligned}$$

where $L(x_0, x_1, \ldots, x_q) = C\|x_0| + \Lambda(x_1,\ldots,x_q)|^4$ for a constant $C$. Given any $\varepsilon > 0$ choose $M > 0$ such that $\int_{\mathbb{R}^{q+1}} L(x_0, x_1, \ldots, x_q) 1_{\{L(x_0,x_1,\ldots,x_q)>M\}} \mu_{q+1}(dx_0, dx_1, \ldots, dx_q) < \varepsilon$. Then we can write

$$\sup_{h\in H}|\mathcal{B}_n(h) - \widetilde{\mathcal{B}}_n(h)| \leq M\|\theta_n^{LS} - \theta_0\| + \frac{1}{n}\sum_{i=1}^n L(X_i, \mathbf{X}_i) 1_{\{L(X_i,\mathbf{X}_i)>M\}} \operatorname{diam}(\Theta).$$

The second term of the right-hand side converges to a positive constant which is smaller than $\varepsilon \cdot \operatorname{diam}(\Theta)$. Since the choice of $\varepsilon > 0$ is arbitrary, we have $\sup_{h\in H}|\mathcal{B}_n(h) - \widetilde{\mathcal{B}}_n(h)| = o_{P^*}(1)$.

Now, we can write $\widetilde{\mathcal{B}}_n(h) = T_{1,n} + T_{2,n}(h) + T_{3,n}(h)$, where

$$T_{1,n} = \frac{1}{n}\sum_{i=1}^n ||X_i - S(\mathbf{X}_i;\theta_0)|^2 - \sigma^2(\mathbf{X}_i;h_0)|^2,$$

$$T_{2,n}(h) = \frac{2}{n}\sum_{i=1}^n \sigma^2(\mathbf{X}_i, h_0)(w_i^2 - 1)(\sigma^2(\mathbf{X}_i, h_0) - \sigma^2(\mathbf{X}_i;h)),$$

$$T_{3,n}(h) = \frac{1}{n}\sum_{i=1}^n |\sigma^2(\mathbf{X}_i;h_0) - \sigma^2(\mathbf{X}_i;h)|^2.$$

The term $T_{1,n}$ converges in probability to a constant $C_1$ by assumption. By Jain–Marcus' CLT for martingale difference arrays, it is easy to show that $\sup_{h\in H}|T_{2,n}(h)|=o_{P^*}(1)$ (here, we use the technical assumption that $E[w_i^4|\mathcal{F}_{i-1}]$ is bounded). Finally, by Theorem 3.1, it holds that $\sup_{h\in H}|T_{3,n}(h)-T_3(h)|=o_{P^*}(1)$ where

$$T_3(h) = \int_{\mathbb{R}^q} |\sigma^2(\mathbf{x};h_0) - \sigma^2(\mathbf{x};h)|^2 \mu_q(d\mathbf{x}).$$

Consequently, we have $\sup_{h\in H}|\mathcal{B}_n(h)-(C_1+T_3(h))|=o_{P^*}(1)$. Therefore, the claim of the second step follows form van der Vaart and Wellner's (1996) consistency theorem. □

**Acknowledgments.** I thank the Associate Editor and the referees for their careful reading and the advice to present a sharper version of the moment conditions.

Institute of Statistical Mathematics
4-6-7 Minami-Azabu, Minato-ku
Tokyo 106-8569
Japan
E-mail: nisiyama@ism.ac.jp
URL: http://www.ism.ac.jp/~nisiyama/